# INTRODUCTION

# TECHNOLOGIES AND RESOURCES IN MATHEMATICAL EDUCATION

Ghislaine Gueudet (CREAD, IUFM Bretagne UBO, France)
Rosa Maria Bottino (ITD, CNR Genova, Italy)
Giampaolo Chiappini (ITD, CNR Genova, Italy)
Stephen Hegedus (Kaput Center, University of Massachusetts, USA)
Hans-Georg Weigand (University of Wuerzburg, Germany)

**INTRODUCTION**

Technologies in mathematical education has been a theme present at CERME from the first edition. The available technologies have evolved a lot during these years. At CERME 5 conference, the conclusions of the technology Working Group (Kynigos et al. 2007), as well as Artigue's and Ruthven's interventions (Artigue 2007, Ruthven 2007), signal perspective evolutions towards more comprehensive studies, in several respects. Drawing on these previous works, CERME 6 WG7 intended to go further in the directions they have indicated.

An important issue, accounting for the introduction of the word "resources" in the name of a group which was previously called "tools and technologies in mathematical didactics", is the need for considering technologies within a range of resources available for the students, the teachers, teacher's trainers etc. These agents can draw on software, computers, interactive whiteboards, online resources, but also on more traditional geometry tools, textbooks, etc. Various kinds of digital material are now extensively used, and they can be viewed as belonging to a wider set of curriculum material (Remillard 2005) and teaching resources (Adler 2000). The papers in WG7 concern different kinds of resources, still with a specific focus on digital material. Another specific focus of WG7 is on theoretical approaches. Design issues need to focus on integration and impact, especially in the use of innovative technology. This entails the development of approaches framing research on fidelity, efficacy, and effective integration (Hegedus & Lesh, 2008). These approaches have been discussed in the group, and several issues linked with the articulation between research and development have been raised, as it is presented below.

The work in the group was organized into three parallel sessions, corresponding to three themes summarized below; specific slots were devoted to the presentation of the work done within three projects co-funded by the European Community, whose participants are represented in WG7: the Telma European Research Team, the Remath project, and the Intergeo project. The whole group was nevertheless gathered for the first session, with two important activities: the identification of questions considered important for the group's work by the participants (figure 1); and the "plenary" address of Jean-Baptiste Lagrange on the results of the ICMI17 study (Hoyles & Lagrange, to appear). The following trends in research and salient elements presented by Jean-Baptiste Lagrange were extensively present in the group's discussions:

− The integration and synthesis of previously fragmented theories and the development of broad approaches;
− The consideration of the design of tools and curricula as a major issue for mathematical education;
− The development of teacher-oriented research studies with a specific consideration on methodological issues such as the consideration of "ordinary" teachers by researchers

> 1. The word "dynamic" permeates mathematics teaching and learning activities which involve technology. Where and when is it appropriate?
> 2. How is it possible to restructure the maths curriculum to take advantage of new technologies to generate mathematical thinking?
> 3. How can we assess the applicability and the effectiveness of current theories?
> 4. Do we need specific theoretical tools or approaches to study the different ways in which teaching can be carried out using technologies?
> 5. What kind of professional development could support pre/in-service teachers to integrate new technologies in their classroom practice?
> 6. How can "old" and "new" resources interact each other? For example, how is it possible to incorporate e-technologies in textbooks?
> 7. If we take seriously into account a semiotic perspective considering the evolution of ICT, what new is offered in terms of creative power of semiotic means?
> 8. Can the European projects presented at WG7 contribute to create a "general theory" of teaching and learning with ICT which can be useful in different European countries?

Figure 1. Examples of questions raised by WG7 participants

## THEMES AND PROJECTS IN WG7

*Design, articulation of design and use*

The work within the « design » theme extensively dealt with the link between design and learning. In particular, the question on the way in which mathematical knowledge can be modified according to different environments was raised. Changes induced by vizualization (Ladel & Kortenkamp, Lois & Milevicich, Kortenkamp & Rolka), virtual reality and simulations (Dana-Picard et al., Bessot & Laborde) opportunities were discussed. Such changes are linked with the specific software tools considered and related features and their analysis is central for the design process. Tools can modify knowledge and learning as well, as it clearly happens for outdoor activities (Nilsson et al.). A broad view on the technology involved, but also on the appropriate associated mathematical tasks is essential (Buteau & Muller, Diakoumopoulos). Amongst the possibilities offered by technology and likely to affect learning and mathematical knowledge, the collective dimensions deserve a specific attention. Software connectivity features, enabling students' collaboration, modify their participation in the mathematical work (Geraniou et al., Hegedus et al.). Questioning collective dimensions also includes a focus on the link between designers and users, and the possible interventions of the users within the design process. This is one of the aspects tackled by the Inter2geo project that dealt with the interoperability of digital geometry systems in Europe (Kortenkamp et al., Trgalova et al.). Beyond the interoperability and indexation issues, this project produced resources for teachers. These resources were tested by users, and users feed-back was included in the design process, with a quality objective. The question of quality, that is the way to assess and validate design and use of technology, was central in the group's discussions. There is still a need for methods to evaluate efficacy of a given technology as far as the learning, engagement and motivation of the students. Even if the projects presented within the group were all grounded in research and linked to given theoretical framework, many questions remain open on the way in which research results can operatively contribute to successful design and use of educational technology.

*Technologies, tools and students mathematical activity*

Papers presented under this theme concerned a wide range of tools and technologies: online resources, software tools, as well as more traditional tools, such as textbooks. Amongst these technologies, Computer Algebra Systems were considered in several papers (Artigue & Bardini, Buteau et al., Weigand & Bichler). But even with a specific interest on CAS, the research works presented in WG7 consider in fact new complex artefacts, articulating CAS and graphing tools in particular, and raised the problem of designing resources to scaffold the use of these artefacts.

Another direction in which the papers showed richness and variety is the theoretical frameworks they draw upon. On the one hand, general theories of mathematical education were considered: functional thinking (Hoffkamp), situated learning, activity theory (Fernandes et al., Jacinto et al.), theory of didactic situations (Aldon & Durand-Guerrier). These theories were used to enlighten specific aspects of learning with technologies: the idea of functional dependency, the mediation provided by an artefact, the didactic contract. On the other hand, several papers refer to specific theories such as that of instrumental approach (Iranzo & Fortuny, Martigone & Antonini, Rezat). This framework, specially designed to study teaching and learning phenomena involving technology, proposes a genesis perspective on learning with technology. It leads to analyses of learning phenomena in terms of schemes. In WG7, precise classification of schemes and of operational invariants were discussed. Such discussions can contribute to a further progress of the instrumental approach framework.

*Interactions between resources and teachers' professional practice*

The acknowledgment of difficulties linked with the integration of technology in classrooms, identified in previous CERME conferences (Kynigos et al. 2007, Drijvers et al. 2006) was still present in CERME 6 WG7 together with the acknowledgment of the key role played by teachers. The need for investigating teachers' beliefs about technology adoption (Chrysostomou & Mousoulides) was recognized as well as the need for conceptualising systemic innovations of educational systems (Ulm). (Emprin, Cantürk-Günhan & Ozen, Faggiano) discussed the importance of setting up pre-service and in-service teachers' training programs taking into account, for in-service teachers, their pre-established repertoires of resources. An evolution from specific studies of individual teachers' practice to investigations of general integration issues was observed, thus moving a step towards theoretical evolutions. As the matter of fact, for example, teacher's use of ICT was examined with a semiotic mediation perspective in (Maracci & Mariotti), while some authors addressed the development of the instrumental approach to study the role of the teacher, drawing on the notion of instrumental orchestration (Trouche 2004), and introducing the consideration of teachers' instrumental genesis (Billington, Bretscher, Drijvers et al.). The acknowledgment of the variety of resources involved in the teacher's activity as well as the need to take into account the whole classroom context, led some authors to develop holistic approaches, such as a documentary approach to didactics (Gueudet & Trouche) and key structuring features of technology integration in the classroom practice (Ruthven).

Delicate methodological issues are attached to the implementation of these theoretical developments, in particular to the question of the "ordinary teacher", which remains open.

*TELMA/Remath projects*

The topic of the articulation of different theoretical frames is central in the TELMA and Remath European projects. In an effort for overcoming the national specificities, these projects developed a cross-experimentation methodology: The key idea around which this methodology was built was the design and the implementation by each team involved in one of these projects of experiments, carried out in real classroom settings, making use of an ICT-based tool developed by another team (Bottino et al., 2009). They also designed meta-tools, in particular scenarios for researchers and for teachers, and proposed developing an integrated theoretical framework (Bottino & Cerulli; Chiappini & Pedemonte; Maracci et al., Markopoulos et al., Moustaki et al., Trgalova & Chaachoua). In fact the three themes of WG7 are present within these large projects, which opened promising methodological and theoretical directions for research.

## CONCLUSION

What do we retain from the work in WG7? The three themes proposed for the contributions, oriented towards design, students and teachers were from the beginning presented as articulated. The design loops integrate more and more the users, students or teachers. The interactions between students and teachers in class are a focus of attention for the researchers. The articulations between different kinds of resources were also extensively discussed, confirming the need for a broad point of view on resources. Research presented in WG7 is focused on technology, but technology does not mean here a precise delimited tool; it includes meta-tools and complex sets of resources. Reflecting on this evolving meaning of technology can be a direction for the work in future CERME conferences.